\documentclass[a4paper, 12pt]{article}

\usepackage[sort&compress]{natbib}
\bibpunct{(}{)}{;}{a}{}{,} 

\usepackage{amsthm, amsmath, amssymb, mathrsfs, multirow, url, subfigure}
\usepackage{graphicx} 
\usepackage{ifthen} 
\usepackage{amsfonts}
\usepackage[usenames]{color}
\usepackage{fullpage}

\RequirePackage[OT1]{fontenc}
\RequirePackage[colorlinks]{hyperref}
\RequirePackage{hypernat}

\theoremstyle{plain}
\newtheorem{theorem}{Theorem}

\makeatletter
\newcommand*{\defeq}{\mathrel{\rlap{%
                     \raisebox{0.3ex}{$\m@th\cdot$}}%
                     \raisebox{-0.3ex}{$\m@th\cdot$}}%
                     =}
\makeatother
\usepackage{times}
\usepackage{bm}
\usepackage{epstopdf}

\newcommand{\sse}{\text{\sc sse}}

\newcommand{\eps}{\varepsilon}
\newcommand{\RR}{\mathbb{R}}

\title{On overfitting and post-selection uncertainty assessments}

\author{Liang Hong \\
Department of Mathematics, Robert Morris University \\
\url{hong@rmu.edu} \\
\mbox{} \\
Todd A. Kuffner\\
Department of Mathematics, Washington University in St. Louis \\
\url{kuffner@wustl.edu} \\
\mbox{} \\
Ryan Martin \\
Department of Statistics, North Carolina State University \\
\url{rgmarti3@ncsu.edu}
}
\date{}

\begin{document}

\maketitle

\begin{abstract}
In a regression context, when the relevant subset of explanatory variables is uncertain, it is common to use a data-driven model selection procedure.  Classical linear model theory, applied naively to the selected sub-model, may not be valid because it ignores the selected sub-model's dependence on the data.  We provide an explanation of this phenomenon, in terms of overfitting, for a class of model selection criteria.  
\smallskip

\end{abstract}


\section{Introduction}
\label{S:intro}

Consider the classical multiple linear regression model 
\begin{equation}
\label{eq:reg}
y = X \beta + \sigma \eps, 
\end{equation}
where $y$ is a $n$-vector of response variables, $X$ is a $n \times p$ matrix of explanatory variables, $\beta$ is $p$-vector of slope coefficients, and $\eps$ is a $n$-vector of independent Gaussian noise.  We assume that $p < n$ and that $y$ and the columns of $X$ are centered so that the intercept term can be ignored. Formally, the model corresponds to the family of distributions \eqref{eq:reg} indexed by $\theta=(\beta,\sigma)$ in $\Theta = \RR^p \times (0,\infty)$.  

In practice, there is often uncertainty about the set of explanatory variables to be included.  In such cases, it is common to express the parameter $\theta$ as $(S, \beta_S, \sigma_S)$, where $S \subseteq \{1,\ldots,p\}$ represents a subset of the explanatory variables, $\beta_S \in \RR^{|S|}$ represents the coefficients corresponding to the specific set $S$, and $\sigma_S > 0$.  This amounts to decomposing the full parameter space $\Theta$ as $\Theta = \bigcup_S \Theta(S)$, where $\Theta(S) = \RR^{|S|} \times (0,\infty)$.  Then the model selection problem boils down to choosing a satisfactory sub-model $\Theta(S)$ or, equivalently, a subset $S$.  Standard tools for carrying out this selection step include the Akaike information criterion, \textsc{aic} \citep{Akaike:1973}, and the Bayesian information criterion, \textsc{bic} \citep{Schwarz:1978}.  These are designed to produce models that suitably balance parsimony and fit.

After a subset $S \subseteq \{1,\ldots,p\}$ of explanatory variables is selected, a secondary goal is to make inference on $S$-specific model parameters $(\beta_S, \sigma_S)$, or functions thereof, and/or predict future values of the response.  A naive approach, recommended in textbooks and commonly used by practitioners, is to replace $X$ in \eqref{eq:reg} with $X_S$, the matrix with only the columns corresponding to $S$, and apply classical normal linear model theory.  For example, for a given $x \in \RR^p$, the classical $100(1-\alpha)$\% confidence interval
\begin{equation}
\label{eq:ci}
C_\alpha(x; S) = x_S^\top \hat\beta_S \pm t_{n-|S|-1}(\alpha/2) \hat \sigma_S \{x_S^\top (X_S^\top X_S)^{-1} x_S\}^{1/2}, 
\end{equation}
can be used for inference on the mean response at the given $x$.  However, as is now well-known \citep{BBBZZ:2013}, the properties that these classical procedures enjoy for a fixed/true $S$ may not hold for a data-dependent choice, $\hat S$.  For example, $C_\alpha(x; \hat S)$ may not have coverage probability equal to $1-\alpha$.    

This note provides an explanation of this lack-of-validity phenomenon by showing that, when the sub-model is selected according to information criteria such as \textsc{aic} and \textsc{bic}, if the selected sub-model overfits, i.e., contains a superset of the explanatory variables in the true model, then the corresponding estimate of the error variance will be smaller than that for the true model. This explains the empirical findings in \citet{HKM:2017}, where prediction intervals based on the sub-model minimizing \textsc{aic} tend to be too short compared to those based on the true model and, consequently, they tend to undercover; see Section~\ref{S:example}.  Moreover, our Theorem~\ref{thm:overfit} together with the dilation phenomenon described in \citet{Efron:2003}, explains why bootstrap may not correct the selection effect for methods that tend to overfit.  

\section{Result}
\label{S:result}

For a given sub-model $\Theta(S)$, corresponding to a subset $S \subseteq \{1,\ldots,p\}$, let $(\hat\beta_S, \hat\sigma_S)$ denote the least squares estimators of the $\Theta(S)$-specific parameters $(\beta_S, \sigma_S)$.  We consider a selection procedure that chooses the subset $S$ by minimizing the function 
\begin{equation}
\label{eq:gamma}
\gamma_n(S) = n \log \sse(S) + c_n |S|, \quad S \subseteq \{1,\ldots,p\}, 
\end{equation}
where $\sse(S) = \|y - X_S \hat\beta_S\|^2$ is the error sum of squares for sub-model $\Theta(S)$, which is proportional to the corresponding least squares estimator $\hat\sigma_S^2$, $c_n=o(n)$ is a user-specified sequence of constants, and $|S|$ denotes the cardinality of the set $S$.  The \textsc{aic} and \textsc{bic} set $c_n \equiv 2$ and $c_n = \log n$, respectively.  

Suppose that there exists a subset $S^\star$ corresponding to the truly non-zero regression coefficients, i.e., $\beta_i \neq 0$ for $i \in S^\star$ and $\beta_i = 0$ for $i \not\in S^\star$.  We write $(\hat\beta_{S^\star}$, $\hat\sigma_{S^\star})$ for the oracle estimators, those based on knowledge of the true sub-model $\Theta(S^\star)$.  Of course, if $\hat S$ is the subset chosen by minimizing $\gamma_n$ in \eqref{eq:gamma}, then $\gamma_n(\hat S) \leq \gamma_n(S^\star)$ or, equivalently, 
\begin{equation}
\label{eq:selection}
n \log \sse(\hat S) + c_n|\hat S| \leq n \log \sse(S^\star) + c_n|S^\star| ;
\end{equation}
if $\hat S \neq S^\star$, then the inequality in \eqref{eq:selection} would be strict. 

For the purpose of inference or prediction, it is common to naively use the classical normal linear model theory, based on the selected subset $\hat S$, to derive uncertainty assessments.  However, using the data to select $\hat S$ introduces bias, violating the assumptions of that classical theory, and thereby invalidating the conclusions.  The next result provides an explanation for this general phenomenon in cases where the selected sub-model $\Theta(\hat S)$ overfits in the sense that $\hat S \supset S^\star$.  In such cases, we find that $\hat\sigma_{\hat S}$ is smaller than the oracle estimator $\hat\sigma_{S^\star}$.  Since the error variance estimate is involved in all uncertainty assessment calculations, and since it is common for selection methods to overfit, especially those based on \textsc{aic} \citep{HurvichTsai:1989}, this systematic under-estimation explains the general lack of validity of the classical inferential tools applied naively in a post-selection context.  

\begin{theorem}
\label{thm:overfit}
Suppose $\hat S \supset S^\star$.  If
\begin{equation}
\label{eq:condition}
1 - \exp(-a_n D_n) > D_n, 
\end{equation}
where $a_n = (c_n / n)(n - |S^\star|-1)$ and $D_n = (|\hat S| - |S^\star|)/(n-|S^\star|-1)$, then $\hat\sigma_{\hat S} < \hat\sigma_{S^\star}$. 
\end{theorem}

To gain some intuition about the condition \eqref{eq:condition}, first note that $a_n D_n$ will tend to be small.  In particular, a very conservative bound is $a_n D_n \leq c_n p / n$, which is small for moderate $c_n$ and $n \gg p$.  Next, since $x \mapsto 1-\exp(-ax)$ is convex for $x > 0$ and $a > 0$, we have $1-\exp(-a_n D_n) > a_n D_n$ for all $D_n$ in an interval $(0, d)$, where $d=d(a_n) \in [0, 1)$.  So, to meet \eqref{eq:condition} we need $a_n > 1$ and, again, we have a conservative bound $a_n \geq c_n (n-p-1)/n$, which itself is greater than 1 for $n \gg p$ and $c_n$ not too small.  In particular, if $n \gg p$ and $c_n \equiv 2$ as in the \textsc{aic}, then \eqref{eq:condition} holds.  

\begin{proof}[of Theorem~\ref{thm:overfit}]
Start by writing $\sse(\hat S)$ in terms of $\sse(S^\star)$.  Let $X_{\hat S}$ and $X_{S^\star}$ denote the sub-matrices corresponding to the indicated subsets, and write $P_{\hat S}$ and $P_{S^\star}$ for the respective projections onto their column spaces.  Then Pythagoras' theorem implies that
\begin{equation*}
\sse(\hat S) 
= \sse(S^\star) + Y^\top (P_{S^\star} - P_{\hat S}) Y = (1-r_n) \sse(S^\star),
\end{equation*}
where 
\[ 
r_n = r_n(S^\star, \hat S) = \frac{|\hat S|-|S^\star|}{n - |\hat S|} F_n(S^\star, \hat S), 
\]
and $F_n(S^\star, \hat S)$ is the usual F-statistic for testing the larger $\Theta(\hat S)$ against the smaller $\Theta(S^\star)$.  Consequently, we choose $\hat S$ over the strictly smaller $S^\star$, according to \eqref{eq:selection}, if and only if $r_n > 1 - \exp(-a_n D_n)$.  

Then the above connection between $\sse(\hat S)$ and $\sse(S^\star)$ immediately gives a comparison between the corresponding variance estimates:
\[ \hat\sigma_{\hat S}^2 = \frac{\sse(\hat S)}{n-|\hat S|-1} = \frac{(1 - r_n) \sse(S^\star)}{n-|\hat S|-1} = \frac{n-|S^\star|-1}{n-|\hat S|-1} (1-r_n) \hat\sigma_{S^\star}^2. \]
As above, we find that $\hat \sigma_{\hat S} < \hat\sigma_{S^\star}$ if and only if $r_n > D_n$.  By condition \eqref{eq:condition}, it follows that the lower bound on $r_n$ derived from over-fitting is greater than that derived from the under-estimation.  Therefore, over-fitting implies under-estimation, proving the claim.  
\end{proof}

\section{Illustration}
\label{S:example}

\begin{figure}[t]
\begin{center}
\subfigure[Scatterplot of $\hat\sigma_{\hat S}$ versus $\hat\sigma_{S^\star}$]{\scalebox{0.45}{\includegraphics{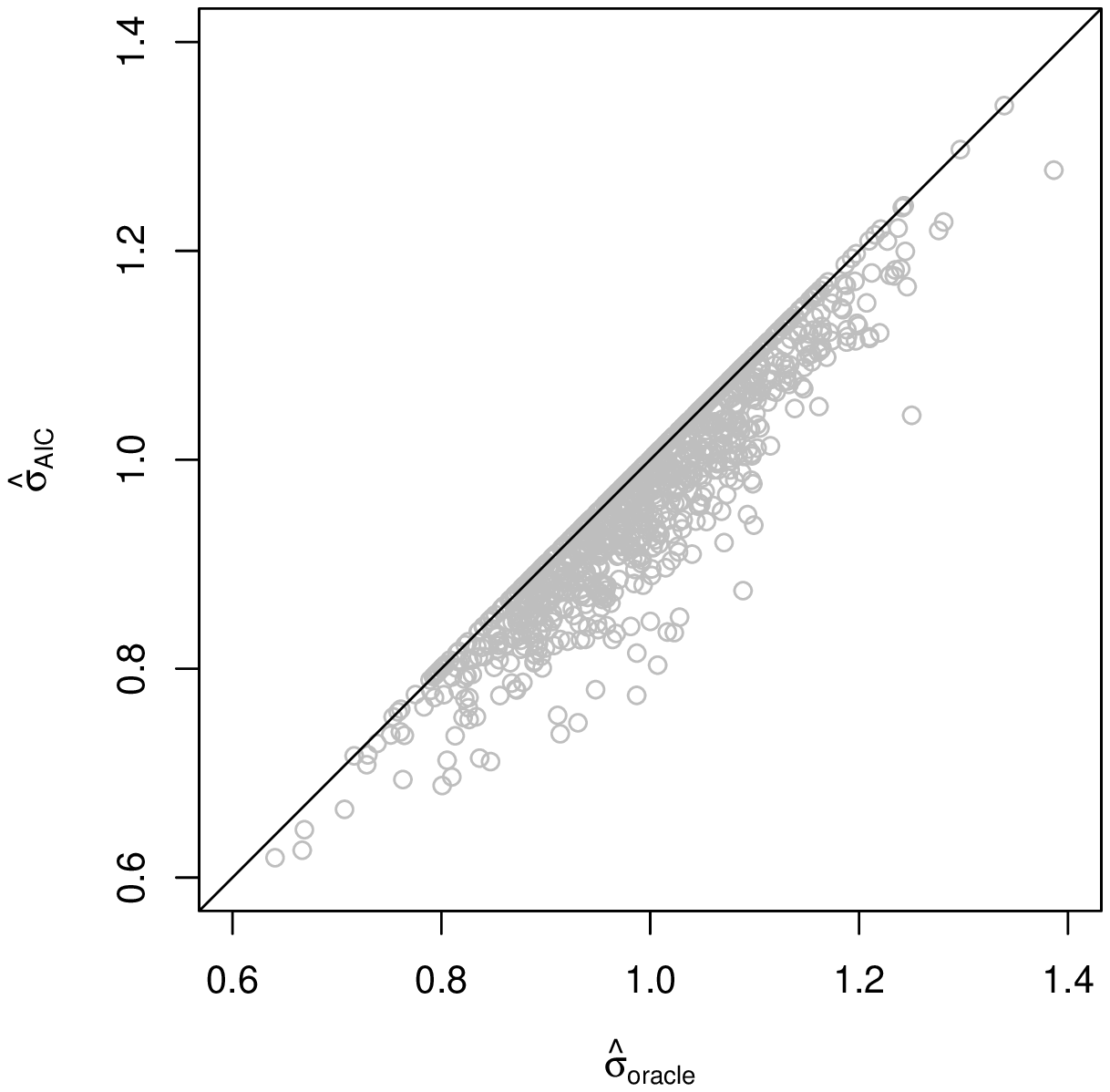}}}
\subfigure[Histogram of $\hat\sigma_{S^\star}/\hat\sigma_{\hat S}$, over-fit only]{\scalebox{0.45}{\includegraphics{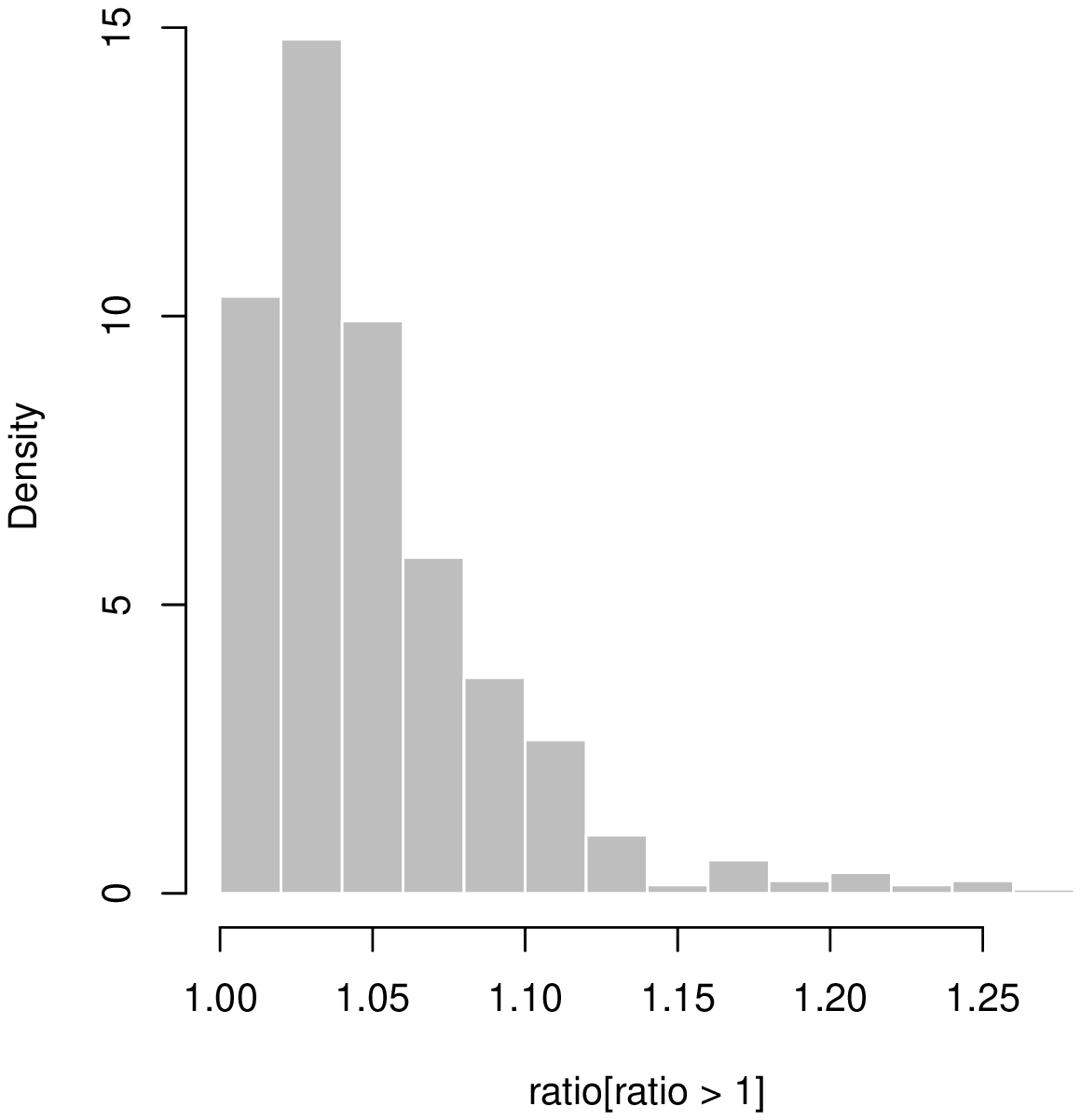}}}
\caption{Plots from the simulations described in Section~\ref{S:example}.}
\label{fig:scatter}
\end{center}
\end{figure}

Consider the model \eqref{eq:reg}, with $n=50$ and $p=10$, and variance $\sigma^2=1$.  Set $S^\star=\{1,2,3\}$, with corresponding coefficients $\beta_1^\star=1$, $\beta_2^\star=2$, and $\beta_3^\star=3$.  The rows of the $X$ matrix are independent, $p$-variate normal, with mean zero, AR(1) dependence structure, and one-step correlation $\rho=0.5$.  We simulated 1000 data sets and, for each, evaluated $\hat\sigma_{\hat S}$ and $\hat\sigma_{S^\star}$, where $\hat S$ is chosen based on the \textsc{aic}.  The scatterplot shown in Figure~\ref{fig:scatter}(a) demonstrates the systematic under-estimation based on the AIC-selected sub-model, as predicted by Theorem~\ref{thm:overfit}.  In all 1000 cases, we have $\hat S \supseteq S^\star$, and those on the diagonal line correspond to $\hat S = S^\star$.  To further illustrate the difference between the estimates, Figure~\ref{fig:scatter}(b) plots a histogram of the ratio $\hat\sigma_{S^\star}/\hat\sigma_{\hat S}$, only for the strict over-fit cases.  In particular, the mean from this histogram is 1.06.  

While the relative difference between the two estimates does not seem remarkable, even this small of a difference can impact the quality of inference.  For example, consider using the confidence interval \eqref{eq:ci} for inference on the mean response at a particular setting $x$ of the explanatory variables; here, $x$ is an independent sample from the distribution that generated the rows of $X$.  The oracle 95\% confidence interval $C_{0.05}(x; S^\star)$ has coverage exactly equal to $0.95$ but, in the 1000 simulations above, the coverage probability of $C_\alpha(x; \hat S)$ is roughly 0.86.  It happens that the $\hat S$-based intervals tend to be shorter than the oracle, suggesting that valid post-selection inference on the mean response requires $\hat\sigma_{\hat S}$ to be strictly larger than $\hat\sigma_{S^\star}$, which is impossible given Theorem~\ref{thm:overfit} and the \textsc{aic}'s tendency to over-fit.

\section*{Acknowledgement}
The authors are grateful to the Editor and two referees whose comments greatly enhanced the clarity of our presentation. Kuffner was supported by the National Science Foundation, U.S.A.

\end{document}